\documentclass{article}
\usepackage[utf8]{inputenc}
\usepackage[english]{babel}
\usepackage{amsthm}
\usepackage{graphicx}
\usepackage{url}
\usepackage{hyperref}
\usepackage{amsmath}

\begin{document}

\title{Asymptotic formulas for harmonic series in terms of a non-trivial zero on the critical line}

\author{Artur Kawalec}

\date{}
\maketitle

\begin{abstract}
In this article, we develop two types of asymptotic formulas for harmonic series in terms of single non-trivial zeros of the Riemann zeta function on the critical line. The series is obtained by evaluating the complex magnitude of an alternating and non-alternating series representation of the Riemann zeta function. Consequently, if the asymptotic limit of the harmonic series is known, then we obtain the Euler-Mascheroni constant with $\log(k)$. We further numerically compute these series for different non-trivial zeros. We also investigate a recursive formula for non-trivial zeros.
\end{abstract}

\section{Introduction}

The harmonic series is defined by a $k$th partial sum of reciprocal positive integers as

\begin{equation}\label{eq:20}
H_k = \sum_{n=1}^{k}\frac{1}{n}.
\end{equation}
It can be shown by Euler-Maclaurin summation formula that the asymptotic expansion of the harmonic series is

\begin{equation}\label{eq:20}
H_k = \gamma+\log(k)+\frac{1}{2k}-\sum_{n=1}^{\infty}\frac{B_{2n}}{2n}\frac{1}{k^{2n}},
\end{equation}
where $B_n$ are Bernoulli numbers and $\gamma$ is the Euler-Mascheroni constant, it can then be extracted by the limit as

\begin{equation}\label{eq:20}
\gamma = \lim_{k\to \infty}\left(H_k-\log(k)\right) = 0.577215664901533\dots,
\end{equation}
which is generally slowly convergent, but one can use more terms of (2) to accelerate the convergence [2]. Also, the Riemann zeta function is traditionally defined by the infinite series

\begin{equation}\label{eq:20}
\zeta(s)= \sum_{n=1}^{\infty}\frac{1}{n^s},
\end{equation}
which is absolutely convergent for $\Re(s)>1$. However, if one substitutes $s=1$, then the structure becomes that of the asymptotic harmonic series (1) as $k\to \infty$. The Riemann zeta function has many different representations with different domains of convergence, see [1][3]. One such representation is the Laurent series which extends analytic continuation of the Riemann zeta function to the whole complex plane
\begin{equation}\label{eq:20}
\zeta(s)=\frac{1}{s-1}+\gamma_0+\sum_{n=1}^{\infty}\frac{(-1)^n\gamma_n(s-1)^n}{n!}
\end{equation}
having only a simple pole at $s=1$. The $\gamma_n$ are the Stieltjes constants defined by

\begin{equation}\label{eq:20}
\gamma_n = \lim_{m\to\infty}\left(\sum_{k=1}^{m}\frac{(\log(k))^n}{k}-\frac{(\log(m))^{n+1}}{n+1}\right),
\end{equation}
where $\gamma_0=\gamma$ is just a special case. As a reference we give a few higher order Stieltjes constants: $\gamma_1=-0.072815...$ and $\gamma_2=-0.009690...$. From (5) one can also extract $\gamma$ near the pole as

\begin{equation}\label{eq:20}
\gamma = \lim_{s\to 1}\left(\zeta(s)-\frac{1}{s-1}\right).
\end{equation}

The zeros of the Riemann zeta function are of great importance since they are related to the distribution of prime numbers, and they come in two types. The first type is the trivial zeros, which occur at negative even integers $s=-2$, $-4$, $-6$, and so on. The second type is the non-trivial zeros, which are complex numbers, and are constrained to lie in a critical strip $0<\Re(s)<1$. We denote a $q$th non-trivial zero as $\rho_q=\sigma_q+it_q$. The Riemann Hypothesis asserts that all non-trivial zeros lie on the critical line at $\sigma=1/2$, which to date is still an unsolved problem; however, all zeros found so far do indeed lie on the critical line. Also, Hardy proved that there is an infinite number of zeros on the critical line, and we distinguish these zeros by $\rho^{'}_q=1/2+it^{'}_q$. The first few non-trivial zeros on the critical line have imaginary components $t^{'}_1 = 14.13472514...$, $t^{'}_2 = 21.02203964...$, $t^{'}_3 = 25.01085758...$ (shown to six decimal places) which are found numerically.

Another common form of the Riemann zeta function is the alternating series

\begin{equation}\label{eq:20}
\zeta(s) = \frac{1}{1-2^{1-s}}\sum_{n=1}^{\infty} \frac{(-1)^{n+1}}{n^s},
\end{equation}
whose domain of convergence is valid for $\Re(s)>0$, however, is not absolutely convergent. One should also note that this form is not an exact analytical continuation of the Riemann zeta function at $\Re(s)=1$ due to additional zeros of the constant factor. But if we consider the critical strip region, then all non-trivial zeros of the Riemann zeta functions are also zeros of (8). Hence under these conditions, we can extract the asymptotic formula for harmonic series in terms of individual non-trivial zeros on the critical line, and which is independent on the validity of the Riemann Hypothesis.  We further numerically validate these series by computing the Euler-Mascheroni constant for different non-trivial zeros using (3). We apply a similar argument for the non-alternating series (4) to obtain the second result. A related formula can be found in [5], where the $\gamma$ is expressed in terms of individual zeros of a cosine integral.

\section{The asymptotic harmonic series first type}
We consider evaluating the complex magnitude of alternating series representation of the Riemann zeta function (8) for $s=\sigma+it$ on the critical strip $0<\sigma<1$,  which results in a form

\begin{equation}\label{eq:20}
\mid \zeta(s) \mid^2 = C^2 (A^2+B^2),
\end{equation}
where the constants $A$ and $B$ are the real and imaginary parts of the infinite sum term of (8) as

\begin{equation}\label{eq:20}
\begin{aligned}
A =  \sum_{n=1}^{\infty} \frac{(-1)^{n}}{n^\sigma}\cos(t\log(n)) \\
B =  -\sum_{n=1}^{\infty} \frac{(-1)^{n}}{n^\sigma}\sin(t\log(n)),
\end{aligned}
\end{equation}
and $C$ is a constant

\begin{equation}\label{eq:20}
C^2 = \frac{1}{1+2^{2(1-\sigma)}-2^{2-\sigma}\cos(t\log(2))}.
\end{equation}
We note that since $A$ and $B$ are convergent in the critical strip, their squares are also convergent, and so the complex magnitude (9) is convergent. When one expands the sum of squares of $A$ and $B$, we have a double sum formed by the complex magnitude, and using certain trigonometric identity simplifies the sum of cosine and sine products into a more compact form

\begin{equation}\label{eq:20}
\mid \zeta(s) \mid^2 = C^2 \sum_{m=1}^{\infty}\sum_{n=1}^{\infty} \frac{(-1)^{m}(-1)^{n}}{m^\sigma n^\sigma}\cos(t (\log(m)-\log(n))).
\end{equation}
We then observe that due to the symmetry of the double sum we have a separation into diagonal ($m=n$) and off-diagonal ($m\neq n$) sums

\begin{equation}\label{eq:20}
\mid \zeta(s) \mid^2 = C^2 \left[\sum_{n=1}^{\infty}\frac{1}{n^{2\sigma}}+2\sum_{n=1}^{\infty}\sum_{m=n+1}^{\infty}\frac{(-1)^{m}(-1)^{n}}{m^\sigma n^\sigma}\cos(t \log(m/n))\right],
\end{equation}
and it is assumed the index variables $m$ and $n$ are positive integers starting with $n=1$ and satisfying $m>n$. Now we consider solutions to

\begin{equation}\label{eq:20}
\mid \zeta(s) \mid^2 = 0
\end{equation}
on the critical strip, which are the non-trivial zeros of the Riemann zeta function, since the real and imaginary parts are coupled together by the complex magnitude. It is evident from (13) that if $\sigma=1/2$ implies the solutions must satisfy

\begin{equation}\label{eq:20}
\sum_{n=1}^{\infty}\frac{1}{n}+2\sum_{n=1}^{\infty}\sum_{m=n+1}^{\infty}\frac{(-1)^{m}(-1)^{n}}{\sqrt{mn}}\cos(t^{'}_q \log(m/n))=0,
\end{equation}
since $C>0$ for all $t$, and $t^{'}_q$ is the $q$th imaginary component of the non-trivial zero on the critical line. We note that these sums are divergent, but they must cancel since (12) is convergent. Therefore, if we treat these divergent sums in a limiting sense, then we have
\begin{equation}\label{eq:20}
\lim_{k\to \infty}\left(\sum_{n=1}^{k}\frac{1}{n}+2\sum_{n=1}^{k}\sum_{m=n+1}^{k}\frac{(-1)^{m}(-1)^{n}}{\sqrt{mn}}\cos(t^{'}_q \log(m/n))\right)=0.
\end{equation}
As a result, the term on the left is the asymptotic harmonic series $H_k$, and the second term is the new series. Hence, the asymptotic relation is
\begin{equation}\label{eq:20}
H_k\sim 2\sum_{n=1}^{k}\sum_{m=n+1}^{k}\frac{(-1)^{m}(-1)^{n+1}}{\sqrt{mn}}\cos(t^{'}_q \log(m/n))
\end{equation}
as $k\to\infty$, where in this form, we absorbed the negative sign by the alternating sign index $n+1$. It then follows that the Euler-Mascheroni constant is

\begin{equation}\label{eq:20}
\gamma = \lim_{k\to \infty}\left(2\sum_{n=1}^{k}\sum_{m=n+1}^{k}\frac{(-1)^{m}(-1)^{n+1}}{\sqrt{mn}}\cos(t^{'}_q \log(m/n))-\log(k)\right).
\end{equation}
The index $q$ of a $t'_q$ zero is valid from one to infinity. Thus the connection to the harmonic series and individual non-trivial zeros on the critical line implies an infinitude of such formulas, which are independent on the validity of the Riemann Hypothesis. We summarize numerical computation of this formula in Appendix A.

\section{The asymptotic harmonic series second type}
Another similar formula for the asymptotic harmonic series can be obtained as such. The application of the Euler-Maclaurin summation formula to the infinite series (4) yields an expansion
\begin{equation}\label{eq:20}
\zeta(s)=\sum_{n=1}^{k-1}\frac{1}{n^s}-\frac{k^{1-s}}{1-s}+\frac{1}{2}k^{-s}+\frac{B_2}{2}sk^{-s-1}+...,
\end{equation}
for an arbitrary $k$, but when $k\to\infty$, then it suffices to use the first term in the expansion, which actually subtracts the pole and extends the analytical continuation of (4) to a new domain $\Re(s)>0$ as
\begin{equation}\label{eq:20}
\zeta(s)=\lim_{k\to\infty}\sum_{n=1}^{k-1}\frac{1}{n^s}-\frac{k^{1-s}}{1-s}.
\end{equation}
Therefore, similarly as before, the non-trivial zeros are also zeros of (20), and if we consider the real and imaginary parts on the critical line, then we have a set of two asymptotic equations
\begin{equation}\label{eq:20}
\begin{aligned}
\sum_{n=1}^{k-1}\frac{1}{\sqrt{n}}\cos(t^{'}_q\log(n))=\frac{\sqrt{k}}{(1/2)^2+t^{'2}_q}\left(\frac{1}{2}\cos(t^{'}_q\log(k))+t^{'}_q\sin(t^{'}_q\log(k))\right) \\
\sum_{n=1}^{k-1}\frac{1}{\sqrt{n}}\sin(t^{'}_q\log(n))=\frac{\sqrt{k}}{(1/2)^2+t^{'2}_q}\left(\frac{1}{2}\sin(t^{'}_q\log(k))-t^{'}_q\cos(t^{'}_q\log(k))\right),
\end{aligned}
\end{equation}
as $k\to\infty$. Then, by considering solutions to (14), we compute the complex magnitude by squaring and adding both equations to obtain an asymptotic relationship

\begin{equation}\label{eq:20}
\frac{k}{(\frac{1}{2})^2+t^{'2}_q} \sim \sum_{m=1}^{k-1}\sum_{n=1}^{k-1} \frac{1}{\sqrt{mn}}\cos(t'_q(\log(m)-\log(n))),
\end{equation}
where we have the same double sum as before, but now without the alternating sign. Hence, by separating the double sum into the diagonal ($m=n$) and off-diagonal sums ($m \neq n$) in the limit, it follows that the second form of the asymptotic harmonic series is
\begin{equation}\label{eq:20}
H_k\sim \frac{k+1}{(\frac{1}{2})^2+t^{'2}_q}-2\sum_{n=1}^{k}\sum_{m=n+1}^{k}\frac{1}{\sqrt{mn}}\cos(t^{'}_q \log(m/n)),
\end{equation}
where we changed the variable $k\to k+1$, and hence, the second form of Euler-Mascheroni constant is
\begin{equation}\label{eq:20}
\gamma =\lim_{k\to\infty} \left(\frac{k+1}{(\frac{1}{2})^2+t^{'2}_q}-2\sum_{n=1}^{k}\sum_{m=n+1}^{k}\frac{1}{\sqrt{mn}}\cos(t^{'}_q \log(m/n))-\log(k)\right)
\end{equation}
for $q=1,2,3...$. Again, this result implies an infinitude of such formulas. We summarize the numerical computation of this formula in Appendix A.

\section{On the non-trivial zeros}
In the above formula, the non-trivial zero term is isolated, and so we can extract it to obtain
\begin{equation}\label{eq:20}
t^{'2}_q =\lim_{k\to\infty} (k+1)\left(\gamma+\log(k)+2\sum_{n=1}^{k}\sum_{m=n+1}^{k}\frac{1}{\sqrt{mn}}\cos(t^{'}_q \log(m/n))\right)^{-1}-\left(\frac{1}{2}\right)^2.
\end{equation}
In general, if we define a function $f(t)$ as
\begin{equation}\label{eq:20}
f(t) =\lim_{k\to\infty} \left((k+1)\left(\gamma+\log(k)+2\sum_{n=1}^{k}\sum_{m=n+1}^{k}\frac{1}{\sqrt{mn}}\cos(t\log(m/n))\right)^{-1}-\left(1/2\right)^2\right)^{1/2},
\end{equation}
then it has a recursive property that
\begin{equation}\label{eq:20}
f(t^{'}_q) = t^{'}_q
\end{equation}
for all $q=1,2,3\dots $. We verified this function by computing non-trivial zeros recursively and summarized the results in Appendix B for several zeros. We note that the zeros can be extracted as a function of zeros, but like equation (24) is slow to converge. There is also a similar equation that can be developed from (21) if we solve for $t'_{q}$.  Using the cosine representation in (21) we obtain
\begin{equation}\label{eq:20}
g(t) = \lim_{k\to\infty} \cot(t\log(k))\left(\frac{(1/2)^2+t^2}{\sqrt{k}\cos(t\log(k))} \sum_{n=1}^{k}\frac{1}{\sqrt{n}}\cos(t\log(n))-\frac{1}{2}\right)
\end{equation}
with same property
\begin{equation}\label{eq:20}
g(t^{'}_q) = t^{'}_q
\end{equation}
for all $q=1,2,3\dots$. We next perform a numerical experiment to attempt to compute non-trivial zeros from an initial value $y_0$ by recursively applying $g(t)$ an $n$ number of times as
\begin{equation}\label{eq:20}
t^{'}_q \approx g(g(g(y_0))) \dots n,
\end{equation}
which under certain circumstances will converge, but not always. To illustrate this, we set an initial value $y_0$=14.2 and then recursively apply $g(t)$, we have

\begin{equation}\label{eq:20}
\begin{aligned}
y_1 = g(y_0) =  14.1989979879525 \\
y_2 = g(y_1) =  14.1980396069588 \\
y_3 = g(y_2) =  14.1971203252332
\end{aligned}
\end{equation}
for $k=10^7$. After a few iterations, the value changed slowly, but after $1000$ iterations, we reached $t'_{1}\approx 14.1346677444012$, which is accurate to $3$ decimal places. After a few numerical experiments, we note that most of the time, this method will not converge due to being highly sensitive to $k$ and initial value, which causes a change in sign in the iteration. We could not reproduce this result with the $f(t)$ formula above, as it is slowly convergent.

\section{Appendix A}
In this section, we numerically compute Euler-Mascheroni constant by equation (18) and (24) for different parameters in the Matlab software package. We used non-trivial zeros from LMFDB database [4]. We report the following results. In Table $1$, we used the first zero $t'_1$ and $k$ as powers of $10$ from $1$ to $5$, and note that the convergence for equation (18) is approaching to the correct value and getting better for higher $k$, but very slowly. For equation (24), the convergence is much slower, even for high $k$.

\begin{table}[ht]
\caption{Evaluation of $\gamma$ by Equation (18) and (24) for $t_1^{'}$} 
\centering 
\begin{tabular}{c c c} 
\hline\hline 
$k$ & $\gamma$ Eq. (18) & $\gamma$ Eq.(24)\\ [0.5ex] 
\hline 
$10^1$  & 0.588166547527396 & 0.624430642787654 \\
$10^2$  & 0.579707476081083 & 0.583918804120366 \\
$10^3$  & 0.577465694084099 & 0.580132200473009 \\
$10^4$  & 0.577240665308434 & 0.579756829762655 \\
$10^5$  & 0.577218164898902 & 0.579719325600715
\\ [1ex] 
\hline 
\end{tabular}
\label{table:nonlin} 
\end{table}

In Table $2$ we run the first $10$ zeros at $k=10^5$, and in Table $3$ we run zeros as powers of $10$ from $2$ to $5$ at $k=10^5$ and note that for (18) the convergence is remarkably stable and accurate to five decimal places, even for high zeros. However, for (24) the convergence is very slow for lower zeros, compared with (18), but starts improving for high zeros.  The number of off-diagonal elements of these double sums is $k^2/2-k$, so for $k=10^5$, that is almost $5$ billion. We also tried even higher-order zeros, but eventually, these series fail to converge to a first digit, and so more terms are required.

\begin{table}[ht]
\caption{Evaluation of $\gamma$ by Equation (18) and (24) for $k=10^5$}
\centering 
\begin{tabular}{c c c c} 
\hline\hline 
$q$ & $t_q^{'}$ & $\gamma$ Eq.(18) & $\gamma$ Eq.(24)\\ [0.5ex] 
\hline 
$1$  & 14.1347251417347 & 0.577218164898902 & 0.579719325600715 \\
$2$  & 21.0220396387716 & 0.577218164886766 & 0.578350602290223 \\
$3$  & 25.0108575801457 & 0.577218164913269 & 0.578018818257371 \\
$4$  & 30.4248761258595 & 0.577218164961156 & 0.577759833545594 \\
$5$  & 32.9350615877392 & 0.577218164938410 & 0.577680674428317 \\
$6$  & 37.5861781588257 & 0.577218164922645 & 0.577573695815113 \\
$7$  & 40.9187190121475 & 0.577218164911756 & 0.577518411665233 \\
$8$  & 43.3270732809150 & 0.577218164859838 & 0.577486145253191 \\
$9$  & 48.0051508811672 & 0.577218164927322 & 0.577436775219061 \\
$10$ & 49.7738324776723 & 0.577218164882106 & 0.577421632805789
\\ [1ex] 
\hline 
\end{tabular}
\label{table:nonlin} 
\end{table}

\begin{table}[ht]
\caption{Evaluation of $\gamma$ by Equation (18) and (24) for $k=10^5$}
\centering 
\begin{tabular}{c c c c} 
\hline\hline 
$q$ & $t_q^{'}$ & $\gamma$ Eq.(18) & $\gamma$ Eq.(24)\\ [0.5ex] 
\hline 
$10^2$  & 236.52422966581 & 0.577218164909381 & 0.577228769071846 \\
$10^3$  &1419.42248094599 & 0.577218164787256 & 0.577220079724956 \\
$10^4$  &  9877.78265400550 & 0.577218158790689 & 0.577219836266831 \\
$10^5$  & 74920.827498994 & 0.577217778781408 & 0.577219808522806

\\ [1ex] 
\hline 
\end{tabular}
\label{table:nonlin} 
\end{table}

\section{Appendix B}
In this section, we validate equation (26) for non-trivial zeros. Similarly as before, in Table $4$, we compute the first zero $t'_1$ and use $k$ as powers of $10$ from $1$ to $5$. We note that convergence is approaching to the true value. In Table $5$ we compute the first $10$ zeros at $k=10^5$, and in Table $6$ we compute zeros as powers of $10$ from $2$ to $5$ at $k=10^5$. We note that the high zeros begin to have more error; as a result,  more terms are required to improve convergence.

\begin{table}[ht]
\caption{Evaluation of $t_1^{'}$ by Equation (26) for different $k$} 
\centering 
\begin{tabular}{c c c} 
\hline\hline 
$k$ & $t'_1$ Eq. (26) \\ [0.5ex] 
\hline 
$10^1$  & 30.2497502548065  \\
$10^2$  & 14.2290157794652  \\
$10^3$  & 14.1388506664484  \\
$10^4$  & 14.1350848277514  \\
$10^5$  & 14.1347605815184
\\ [1ex] 
\hline 
\end{tabular}
\label{table:nonlin} 
\end{table}

\begin{table}[ht]
\caption{Evaluation of $t_q^{'}$ by Equation (26) for $k=10^5$}
\centering 
\begin{tabular}{c c c} 
\hline\hline 
$q$ & $t_q^{'}$ & $t'_q$ Eq. (26)\\ [0.5ex] 
\hline 
$1$  & 14.1347251417347 & 14.1347605815185  \\
$2$  & 21.0220396387716 & 21.0220924170205   \\
$3$  & 25.0108575801457 & 25.0109204581271  \\
$4$  & 30.4248761258595 & 30.4249527952168   \\
$5$  & 32.9350615877392 & 32.9351446884539  \\
$6$  & 37.5861781588257 & 37.5862732468959 \\
$7$  & 40.9187190121475 & 40.9188227512751\\
$8$  & 43.3270732809150 & 43.3271833072802\\
$9$  & 48.0051508811672 & 48.0052732114747   \\
$10$ & 49.7738324776723 & 49.7739594934774

\\ [1ex] 
\hline 
\end{tabular}
\label{table:nonlin} 
\end{table}

\begin{table}[ht]
\caption{Evaluation of $t_q^{'}$ by Equation (26) for $k=10^5$}
\centering 
\begin{tabular}{c c c c} 
\hline\hline 
$q$ & $t_q^{'}$ & $t'_q$ Eq. (26) \\ [0.5ex] 
\hline 
$10^2$  & 236.52422966581 & 236.52509664502  \\
$10^3$  &1419.42248094599 & 1419.48561150748 \\
$10^4$  &  9877.78265400550 & 9897.94562749923 \\
$10^5$  & 74920.827498994 & 85523.0271275466

\\ [1ex] 
\hline 
\end{tabular}
\label{table:nonlin} 
\end{table}

\newpage
\section{Conclusion}
We presented a simple asymptotic formula for harmonic series in terms of individual non-trivial zeros on the critical line using the alternating and non-alternating series representation of the Riemann zeta function, from which $\gamma$ follows using the $\log(k)$. The numerical computation of these series shows that the convergence is very slow. We note that equation (18) for $\gamma$ is converging faster than equation (26). It is also more stable across the non-trivial zeros on the critical line. Also, we investigated recursive formulas for non-trivial zeros, which were verified numerically; however, they are not practical. Perhaps this method could be improved to compute non-trivial zeros.

\section*{Acknowledgement}
I would like to thank Dr. Marek Wolf for providing valuable feedback during the development of this article, as well as performing an independent numerical computation of the presented formulas to high precision.

\texttt{Email: art.kawalec@gmail.com}

\end{document}